\documentclass[10pt,reqno]{article}

\usepackage{amssymb}
\usepackage{amsmath}
\usepackage{graphicx}
\usepackage[ruled]{algorithm2e}

\newcommand{\Id}{I}
\newcommand{\lo}{$\mathcal C^1$-Lohner }
\newcommand{\ho}{$\mathcal C^1$-HO }
\newcommand{\intv}[1]{\left[\underline{#1},\overline{#1}\right]}
\newcommand{\tc}[2]{{#1}^{[#2]}}
\newcommand{\dc}[2]{{#1}^{(#2)}}
\def\qed{{\hfill{\vrule height5pt width3pt depth0pt}\medskip}}

\newtheorem{thm}{Theorem}
\newtheorem{lemma}[thm]{Lemma}

\title{An implicit algorithm for validated enclosures of the solutions to 
variational equations for ODEs}
\author{
    {\Large \bf Irmina Walawska and Daniel Wilczak\footnote{
This research is partially supported by 
the Polish National Science Center under Maestro Grant No. 2014/14/A/ST1/00453. 
       }
}\\\\
    Faculty of Mathematics and Computer Science\\
    Jagiellonian University\\
    {\L}ojasiewicza 6, 30-348 Krak\'ow, Poland\\
    \small \texttt{\{Irmina.Walawska, Daniel.Wilczak\}@ii.uj.edu.pl}
}
\date{\today}

\begin{document}

\maketitle

\begin{abstract}
We propose a new algorithm for computing validated bounds for the solutions to
the first order variational equations associated to ODEs. These validated 
solutions are the kernel of numerics computer-assisted proofs in dynamical 
systems literature.
The method uses a high-order Taylor method as a predictor step and an implicit 
method based on the Hermite-Obreshkov interpolation as a corrector step.
The proposed algorithm is an improvement of the $\mathcal C^1$-Lohner algorithm 
proposed by Zgliczy\'nski and it provides sharper bounds.

As an application of the algorithm, we give a computer-assisted proof of the 
existence of an attractor set in the R\"ossler system, and we show that the 
attractor contains an invariant and uniformly hyperbolic subset on which the 
dynamics is chaotic, that is, conjugated to subshift of finite type with 
positive topological entropy.
\end{abstract}


\section{Introduction.}
The aim of this paper is to provide an algorithm that computes validated enclosures for the solutions to the following set of initial value problems
\begin{eqnarray}\label{eq:vivp}
 \left\{\begin{array}{rcl}
 \dot x(t) &=& f(x(t)),\\
 \dot V(t) &=& Df (x(t))\cdot V(t),\\
 x(0) &\in& [x_0],\\
 V(0) &\in& [V_0],
 \end{array}\right.
\end{eqnarray}
where $f:\mathbb R^n\to \mathbb R^n$ is a smooth function  (usually analytic in the domain) and $[x_0]\subset \mathbb R^n$, $[V_0]\subset \mathbb R^{n^2}$ are sets of initial conditions. In contrast to standard numerical methods, one step of a validated algorithm for (\ref{eq:vivp}) produces sets $[x_1]$ and $[V_1]$ that guarantee to contain $x(h)\in [x_1]$ and $V(h)\in [V_1]$ for all initial conditions $x(0)\in[x_0]$ and $V(0)\in [V_0]$, where $h>0$ is a time step (usually variable) of the method. The computations are performed in interval arithmetics \cite{Moore1966} in order to obtain guaranteed bounds on the expressions we evaluate.

The equation for $V(t)$ is called the variational equation associated with an 
ODE. Solutions $V(t)$ give us an information about sensitivities of trajectories 
with respect to initial conditions. They proved to be useful in finding 
periodic 
solutions, proving their existence and analysis of their stability 
\cite{BarrioRodriguez2014,BarrioRodriguezBlessa2012,Galias2006,GaliasTucker2008,
KapelaSimo2007,KapelaZgliczynski2003}. They are used to estimate invariant 
manifolds of periodic orbits \cite{Capinski2012,CapinskiWasieczko2015}. 
Derivatives with respect to initial conditions are used to prove the existence 
of connecting orbits 
\cite{BarrioMartinezSerranoWilczak2015,SzczelinaZgliczynski2013,
WilczakZgliczynski2009siads} or even (non)uniformly hyperbolic and chaotic 
attractors \cite{Tucker2002,Wilczak2010}. First and higher-order derivatives 
with respect to initial conditions can be used to study some bifurcation 
problems \cite{KokubuWilczakZgliczynski2007,
WilczakZgliczynski2009siads,WilczakZgliczynski2009focm}. This wide spectrum of 
applications is our main motivation for developing an efficient algorithm that 
produces sharp bounds on the solutions to (\ref{eq:vivp}).

In principle, the problem (\ref{eq:vivp}) can be solved by any algorithm capable 
to compute validated solution to IVP for ODEs. There are several available 
algorithms and their implementations --- just to mention a few of them: VNODE-LP 
\cite{NedialkovJacksonCorliss1999,Nedialkov2006,NedialkovJackson1998,
NedialkovJacksonPryce2001}, COSY Infinity 
\cite{BerzMakino1999,MakinoBerz2006,MakinoBerz2009}, CAPD \cite{CAPD}, 
Valencia-IVP \cite{RauhBrillGunter2009}. The above mentioned ODE solvers are 
internally higher-order methods with respect to the initial state, which means 
that they use at least partial information about the derivatives with respect to 
initial conditions to reduce the \emph{wrapping effect}. Therefore, direct 
application of these solvers uses a higher effective dimension (the 
internal dimension of the solver) than the dimension of the phase space. In the 
case of the 
codes VNODE-LP and CAPD, this effective dimension is $(n(n+1))^2$, which 
dramatically decreases the performance of these methods when applied directly to 
the extended system (\ref{eq:vivp}). This key observation motivated developing 
the \lo algorithm \cite{Zgliczynski2002}, which takes into account the block 
structure of (\ref{eq:vivp}) and works in $n^2$ effective dimension. Even in low 
dimensions, it is orders of magnitude faster than direct application of a 
$\mathcal C^0$ solver to the variational equations. However, it does not use 
derivatives of $V(t)$ with respect to other coefficients of $V$ (i.e. second 
order derivatives of the original system) to better control the wrapping effect. 
This is why it usually produces worse estimations than those obtained from 
direct application of a $\mathcal C^0$ solver to the extended system.

In this paper, we propose a new algorithm for computation of 
validated solutions to (\ref{eq:vivp}). Our algorithm consists of two steps. 
First, the high-order Taylor method is used as a predictor step. Then, an 
implicit method based on the Hermite-Obreshkov (HO) formula is used to compute 
tighter bounds for the variational equations. This last step is motivated by the 
very famous and efficient algorithm proposed by Nedialkov and Jackson 
\cite{NedialkovJackson1998} and implemented by Nedialkov in the VNODE-LP package 
\cite{Nedialkov2006}. We name the proposed algorithm \ho because it computes 
bounds for the first order variational equations and it is based on the 
Hermite-Obreshkov interpolation formula.

Our algorithm, by its construction, cannot produce worse estimations than the 
\lo algorithm. Complexity analysis (see Section~\ref{sec:complexity}) shows 
that, in low dimensions, it is slower than the \lo algorithm by the factor 
$9/8$ only. This lack of performance is compensated by a significantly smaller 
truncation error of the method. This allows to take larger time steps when 
computing the trajectories and thus our algorithm appears to be slightly faster 
than the \lo in real applications --- see Section~\ref{sec:applications} for 
the 
case study.

We would like to emphasize that the proposed method can be directly extended to 
the nonautonomous case without increasing the effective dimension of the 
problem. For simplicity in the notation, we consider the autonomous case only. 
In \cite{CAPD}, we provide an implementation of the \ho algorithm for the 
nonautonomous case.

As an application of the proposed algorithm we give a computer-assisted proof 
of the following new result concerning the R\"ossler system \cite{Rossler1976}.
\begin{thm}\label{thm:main}
For the parameter values $a=5.7$ and $b=0.2$ the system 
\begin{equation}\label{eq:rosslerSystem}
\left\{
    \begin{array}{rcl}
        \dot x &=& -y-z\\
        \dot y &=& x+by\\
        \dot z &=& b+z(x-a)
    \end{array}
\right.
\end{equation}
admits a compact, connected invariant set $\mathcal A$ that is an attractor.  
There is an invariant subset $\mathcal H\subset \mathcal A$ on which the 
dynamics is uniformly hyperbolic and chaotic, that is, conjugated to a 
subshift of finite type with positive topological entropy.
\end{thm}
Verification that an ODE is chaotic is not an easy task in general. After 
development of rigorous ODE solvers there appeared numerous 
computer-assisted proofs of the existence of chaos in classical low-dimensional 
systems --- just to mention two pioneering results 
\cite{MischaikowMrozek1995,Zgliczynski1997}. To the 
best of our knowledge there are only two computer-assisted proofs of the 
existence of chaotic and (non)-uniformly hyperbolic attractors for ODEs 
\cite{Tucker2002,Wilczak2010}. These results became possible with development 
of suitable theory and the algorithms capable to integrate variational 
equations. In \cite{Wilczak2010} the \lo algorithm implemented in the CAPD 
library \cite{CAPD} was used.

The paper is organized as follows. In Section~\ref{sec:notation}, we introduce 
the symbols and notation used in the paper. Section~\ref{sec:algorithm} contains 
description of the algorithm and the proof of its correctness. In 
Section~\ref{sec:complexity}, we analyze the complexity of the \ho algorithm and 
we compare it to the complexity of the \lo algorithm. In 
Section~\ref{sec:benchmarks}, we compare the bounds obtained by the \lo and \ho 
algorithms on several examples. In Section~\ref{sec:applications}, we give a 
more detailed statement and proof of Theorem~\ref{thm:main}. We discuss also how 
the computing time depends on the choice of the \lo and \ho algorithms to 
integrate variational equations.
\subsection{Notation.} \label{sec:notation}
By $\Id$ we denote the identity matrix of the dimension clear from the context. 
By $Df$ we denote the derivative of a smooth function $f:\mathbb R^n\to \mathbb 
R^m$. By $D_{x}f$ we denote the partial derivative of $f$ with respect to $x$.

The local flow induced by an ordinary differential equation (ODE) $\dot x(t) = 
f(x(t))$ will be denoted by $\varphi$, i.e.  $\varphi(\cdot,x)=x(\cdot)$ is the 
unique solution passing through $x$ at time zero. We will often identify an 
element $x\in\mathbb R^n$ with the function $x(\cdot)=\varphi(\cdot,x)$. We 
denote $\psi(t,x,V)= D_x\varphi(t,x)\cdot V$. Clearly $\psi(\cdot,x,V)$ is a 
solution to the first-order variational equation associated with an ODE with 
the initial conditions $x$ and $V$.

Let $f: \mathbb{R} \to \mathbb{R}^n$ be a smooth function. By $\dc f i (x)$ we denote the vector of $i$th derivatives of $f$. Normalized derivatives (Taylor coefficients) will be denoted by $\tc f i (x) = \frac{1}{i!}\dc f i(x)$. We apply this notation to derivatives and Taylor coefficients of the flows $\varphi$ and $\psi$ taken with respect to the time variable 
\begin{eqnarray*}
\tc \varphi i(t,x) &:=& (\tc \varphi i(\cdot,x))(t),\\
\tc \psi i(t,x,V) &:=& (\tc \psi i(\cdot,x,V))(t).
\end{eqnarray*}
 
Interval objects will be always denoted in square brackets, for instance $[a]=\intv{a}$ is an interval, and $[v]=([v_1],\ldots,[v_n])$ is an interval vector. Matrices or interval matrices will be denoted by capital letters, for example $[A]$. Vectors and scalars will be always denoted by small letters. We also identify Cartesian product of intervals $[v_1]\times\cdots\times[v_n]$ with a vector of intervals $([v_1],\ldots,[v_n])$. Thus interval vectors can be seen as subsets of $\mathbb R^n$. The same identification will apply to interval matrices.

The midpoint of an interval $[a]=\intv{a}$ will be denoted by $\widehat a= 
(\underline{a}+\overline{a})/2$. The same convention will be used to 
denote the midpoint of an interval vector or an interval matrix; for example 
for 
$[v]=([v_1],\ldots,[v_n])$ we put $\widehat v = (\widehat v_1,\ldots,\widehat 
v_n)$. Sometimes, we will 
denote the midpoint of products of interval objects by $\mathrm{mid}([V][r])$.

Throughout this article, interval vectors or interval matrices marked with tilde 
$[\widetilde y]$, $[\widetilde V]$ will always refer to rough enclosures for the 
solutions to the IVP problem (\ref{eq:vivp}) --- see 
Section~\ref{sec:enclosure} 
for details.
\section{The algorithm.}\label{sec:algorithm}
Consider the initial value problem (\ref{eq:vivp}) and assume that we have 
already proved the existence of the solutions at time $t_k$, and we have 
computed rigorous bounds $[x_k]$ and $[V_k]$ for $\varphi(t_k,[x_0])\subset 
[x_k]$, $\psi(t_k,[x_0],[V_0])\subset [V_k]$, respectively. Let us fix a time 
step $h_k>0$. A rigorous numerical method for (\ref{eq:vivp}) consists usually 
of the following two steps:
\begin{itemize}
 \item computation of a rough enclosure. In this step, the algorithm validates 
that solutions indeed exist over the time interval $[t_k,t_k+h_k]$, and it 
produces sets $[\widetilde y]$ and $[\widetilde V]$, called rough enclosures, 
which satisfy
\begin{eqnarray}
\varphi([0,h_k],[x_k])&\subset& [\widetilde y]\label{eq:c0enc} \quad\textrm{ and 
}\\
\psi([0,h_k],[x_k],\Id)&\subset& [\widetilde V]\label{eq:c1enc}.
\end{eqnarray}
 \item computation of tighter bounds $[x_{k+1}]$, $[V_{k+1}]$ satisfying 
$\varphi(t_k+h_k,[x_0]) \subset [x_{k+1}]$ and $\psi(t_k+h_k,[x_0], [V_0]) \subset [V_{k+1}]$. 
\end{itemize}

In the sequel, we give details of each part of the proposed algorithm.

\subsection{Computation of a rough enclosure.}\label{sec:enclosure}
This section is devoted to describe a method for finding rough enclosures (\ref{eq:c0enc}-\ref{eq:c1enc}). The key observation is that the equation for $V$ in (\ref{eq:vivp}) is linear in $V$, which implies that the following identity holds
$$\psi(t_k+h_k,x_0,V_0) = \psi(h_k,\varphi(t_k,x_0),\Id)\cdot 
\psi(t_k,x_0,V_0)$$
provided all quantities are well defined. This implies that
\begin{equation}\label{eq:V-eval}
\psi(t_k+h_k,[x_0],[V_0]) \subset \psi(h_k,[x_k],\Id)\cdot [V_k].
\end{equation}
Hence, it is sufficient to use $\Id$ as an initial condition for the variational 
equations when computing a rough enclosure $[\widetilde V]$.

One approach to find rough enclosures $[\widetilde y]$ and $[\widetilde V]$ is 
to 
compute them separately. Given a set $[\widetilde y]$ satisfying 
(\ref{eq:c0enc}) 
and computed by any algorithm 
\cite{MrozekZgliczynski2000,NedialkovJacksonPryce2001}, we can try to find an 
interval matrix $[\widetilde 
V]$ such that 
\begin{equation}\label{eq:vfoe}
\Id + [0,h_k]Df([\widetilde y])\cdot [\widetilde V]\subset [\widetilde V].
\end{equation}
If we succeed, then the interval matrix $[\widetilde V]$ satisfies 
(\ref{eq:c1enc}). 
This method is known as the First Order Enclosure (FOE). It has, however, at 
least one significant disadvantage. If we already know that the solutions to 
the 
main equations exist over the time step $h_k$ ($[\widetilde y]$ has been 
computed) 
there is no reason to shorten this time step because solutions to variational 
equation also exist over the same time range. However, this shortening might be 
necessary to fulfill the condition (\ref{eq:vfoe}). To avoid this drawback, 
Zgliczy\'nski proposes \cite{Zgliczynski2002} a method based on logarithmic 
norms that always computes an enclosure $[\widetilde V]$ satisfying 
(\ref{eq:c1enc}) 
for the same time step $h_k$, provided we are able to find an enclosure 
$[\widetilde 
y]$ satisfying (\ref{eq:c0enc}). This type of enclosure is also used in the 
$\mathcal C^r$-Lohner algorithm \cite{WilczakZgliczynski2011} for higher order 
variational equations.

Another strategy is to use the High Order Enclosure (HOE) method 
\cite{CorlissRihm1996,NedialkovJacksonCorliss1999,NedialkovJacksonPryce2001}. 
The authors propose to predict a rough enclosure of the form
\begin{equation}\label{eq:encform}
 [\widetilde y] = \sum_{i=0}^m [0,h_k]^i \tc \varphi i (0,[x_k]) + 
[\varepsilon],
\end{equation}
where $[\varepsilon]$ is an interval vector centered at zero. The inclusion
\begin{equation}\label{eq:hoe}
 [0,h_k]^{m+1}\tc \varphi {m+1} (0,[\widetilde y])\subset [\varepsilon]
\end{equation}
implies that the set $[\widetilde y]$ is indeed a rough enclosure, i.e. it 
satisfies 
(\ref{eq:c0enc}). 
If the inclusion (\ref{eq:hoe}) is not satisfied, then we can always find $\bar 
h_{k}<h_k$ such that (\ref{eq:hoe}) holds with this time step, and thus 
$[\widetilde 
y]$ is a rough enclosure for the time step $\bar h_k$. This strategy is very 
efficient because we do not need to recompute $[\widetilde y]$. Furthermore, 
with 
quite high order $m$ and a reasonable algorithm for time step prediction, we 
usually have $\bar h_k/h_k=\left(\frac{\|\varepsilon\|}{\|\tc \varphi {m+1} 
(0,[\widetilde y])\|}\right)^{1/(m+1)}\approx 1$.

The above method can be used to find simultaneously two enclosures 
$([\widetilde 
y],[\widetilde V])$ for the entire system (\ref{eq:vivp}). We predict 
$[\widetilde y]$ 
as in (\ref{eq:encform}) and 
\begin{equation}\label{eq:vencform}
 [\widetilde V] = \sum_{i=0}^m [0,h_k]^i \tc \psi i (0,[x_k],\Id) + [E].
\end{equation}
Then we check simultaneously (\ref{eq:hoe}) and 
\begin{equation}\label{eq:vhoe}
 [0,h_k]^{m+1}\tc \psi {m+1} (0,[\widetilde y],\Id)[\widetilde V]\subset [E].
\end{equation}
Due to linearity of the equation for variational equations we can consider two 
strategies when (\ref{eq:vhoe}) is not satisfied. We can 
\begin{enumerate}
\item either shorten the time step as we do in (\ref{eq:hoe}) or
\item compute $[\widetilde V^0]$ from the logarithmic norms for the same time 
step 
$h_k$ as proposed in \cite{Zgliczynski2002} and set $[E]=[0,h_k]^{m+1}\tc \psi 
{m+1} (0,[\widetilde y],\Id)[\widetilde V^0]$. Then $[\widetilde V]$ computed as 
in 
(\ref{eq:vencform}) with such $[E]$ satisfies (\ref{eq:c1enc}).
\end{enumerate}
We would like to emphasize that in both cases we do not need to recompute 
the coefficients $\tc \psi 
i 
(0,[\widetilde y],\Id)$ which is very expensive. The first strategy is 
recommended 
when the tolerance per one step is specified which means that there is a 
maximal 
norm of $[E]$ which should not be exceeded. The second strategy applies when 
the 
fixed time step is used (by the user choice or application specific reason).

Our tests show that the $\mathcal C^1$ version of (HOE) gives better results 
than the approach proposed by Zgliczy\'nski, which uses logarithmic norms. Since 
computing a rough enclosure (\ref{eq:c1enc}) is not the main topic of the paper 
we omit details here. 

In what follows we assume that we have a routine that returns three quantities: 
$h_k$, $[\widetilde y]$, $[\widetilde V]$ for which the properties 
(\ref{eq:c0enc}) and (\ref{eq:c1enc}) hold.

\subsection{The predictor step.}
We give a short description of the \lo algorithm \cite{Zgliczynski2002} which will be used as a predictor step in the \ho algorithm. 

\begin{algorithm}[tbp]
 \SetKwInOut{Input}{Input}
 \SetKwInOut{Output}{Output}
 
 \Input{$m$ - natural number (order of the Taylor method)\\
 $h_k$ - positive real number (current time step)\\
 $[x_k]$, $[\widetilde y]$ - interval vectors\\
 $[\widetilde V]$ - interval matrix}
 \Output{$([x_{k+1}^0],[r^0],[V^0],[R^0])\subset\mathbb R^n\times\mathbb 
R^n\times\mathbb R^{n^2}\times\mathbb R^{n^2}$}
  \textbf{Compute}:\\


  $[A]\gets \sum_{i=0}^mh_k^i \tc \psi i (0,[x_k],\Id)$\;
  $y_0\gets\sum_{i=0}^mh_k^i \tc \varphi i (0,\widehat x_k)$\;
  $[y]\gets\sum_{i=0}^mh_k^i \tc \varphi i (0,[x_k])$\;
  $[r^0]\gets [0,h_k]^{m+1}\tc \varphi {m+1}(0,[\widetilde y])$\;
  $[x_{k+1}^0]\gets \left(y_0 + [A]([x_k]-\widehat x_k)\right)\cap [y] + 
[r^0]$\;
  $[R^0]\gets [0,h_k]^{m+1}\tc \psi {m+1}(0,[\widetilde y],\Id)[\widetilde V]$\;
  $[V^0] \gets [A] + [R^0]$\;
 \Return{$([x_{k+1}^0],[r^0],[V^0],[R^0])$}\;
 \caption{Predictor.}\label{alg:predictor}
\end{algorithm}
\begin{lemma}
Assume $h_k$, $[x_k]$, $[\widetilde y]$, $[\widetilde V]$ are such that 
(\ref{eq:c0enc}) 
and (\ref{eq:c1enc}) hold. Then the quantities 
$([x_{k+1}^0],[r^0],[V^0],[R^0])$ 
computed by the Algorithm~\ref{alg:predictor} satisfy
\begin{eqnarray}
 \varphi(h_k,[x_k]) & \subset & [x_{k+1}^0],\label{eq:predictor-property1}\\
 \psi(h_k,[x_k],\Id) & \subset & [V^0],\\
 \mbox{}[0,h_k]^{m+1} \tc \varphi {m+1}(0,[\widetilde y]) & \subset & [r^0] 
\quad
\textrm{ 
and }\\
 \mbox{}[0,h_k]^{m+1} \tc \psi {m+1}(0,[\widetilde y],[\widetilde V]) & \subset 
& 
[R^0]\label{eq:predictor-property2}.
\end{eqnarray}
\end{lemma}
\textbf{Proof:}
The Taylor theorem with Lagrange remainder implies that for all $x_k\in[x_k]$ 
and each component $j=1,\ldots,n$:
\begin{equation}\label{eq:taylor-sum}
\varphi_j(h_k,x_k) = \sum_{i=0}^mh_k^i\tc {\varphi_j} i (0,x_k) + h_k^{m+1}\tc 
{\varphi_j} {m+1}(\tau_j,x_k)
\end{equation}
for some $\tau_j\in(0,h_k)$. By the assumptions 
$\varphi(\tau_j,x_k)\in[\widetilde 
y]$ and by the group property of the flow, we have
$$
  \tc {\varphi_j} {m+1}(\tau_j,x_k) = \tc {\varphi_j} {m+1} 
(0,\varphi(\tau_j,x_k)) \in \tc {\varphi_j} {m+1} (0,[\widetilde y]).
$$
Therefore 
$$
\varphi(h_k,x_k) \in \sum_{i=0}^mh_k^i\tc {\varphi_j} i (0,x_k) + h_k^{m+1}\tc 
{\varphi} {m+1}(0,[\widetilde y]) \subset [y] + [r^0].
$$

Since $[x_k]$ is convex, we can apply the mean value form to the polynomial 
part 
of (\ref{eq:taylor-sum}) and obtain that for $x_k\in[x_k]$ there holds
\begin{equation*}
\sum_{i=0}^mh_k^i\tc \varphi i (0,x_k) \in \sum_{i=0}^mh_k^i\tc \varphi i 
(0,\widehat x_k) + [A](x_k-\widehat x_k)\subset y_0+[A]([x_k]-\widehat 
x_k).
\end{equation*}
Gathering the above together, we obtain
\begin{equation}\label{eq:predictor-eval}
\varphi(h_k,[x_k]) \subset \left(y_0 + [A]([x_k]-\widehat x_k)\right)\cap 
[y] 
+ 
[r^0] = [x_{k+1}^0].
\end{equation}
In a similar way, we deduce that for $x_k\in[x_k]$ and for each component 
$j,c=1,\ldots,n$ there holds
\begin{equation*}
 \tc {\psi_{j,c}} {m+1} (\tau_{j,c},x_k,\Id) \in \tc {\psi_{j,c}} {m+1} 
(0,[\widetilde y],[\widetilde V])
\end{equation*}
for $j,c=1,\ldots,n$, and in consequence
\begin{equation*}
\psi(h_k,[x_k],\Id) \subset [A] + [R^0] = [V^0].
\end{equation*}
\qed
 
\subsection{The corrector step.}
The goal of this section is to set forth a one-step method that refines the 
results obtained from the predictor step and returns tighter rigorous bounds for 
the solution to the ODE and its associated variational equation (\ref{eq:vivp}). 
The method combines the algorithm by Nedialkov and Jackson 
\cite{NedialkovJackson1998} based on the Hermite-Obreshkov interpolation formula 
with the \lo algorithm for variational equations proposed by Zgliczy\'nski 
\cite{Zgliczynski2002}. For reader's convenience, we recall here the key ideas 
of the Hermite-Obreshkov method.

For natural numbers $p,q,i$ such that $i\leq q$, let
\begin{eqnarray*}
    c_i^{q,p} &=& \binom{q}{i}/\binom{p+q}{i}.
\end{eqnarray*}
For a smooth function $u:\mathbb R\to \mathbb R^n$ and real numbers $h,t$ we define
\begin{equation*}
\Psi_{q,p}(h,u,t) = \sum_{i=0}^qc_i^{q,p}h^i \tc u i(t).
\end{equation*}
Using this notation, the Hermite-Obreshkov \cite{Obreschkoff1940} formula reads
\begin{equation}\label{eq:ho-formula}
\Psi_{q,p}(-h,u,h) = \Psi_{p,q}(h,u,0) + (-1)^qc_q^{q,p}h^{p+q+1} R(h,u),
\end{equation}
where
$$R(h,u) = \left(\tc {u_1} {p+q+1}(\tau_1),\ldots,\tc {u_n} 
{p+q+1}(\tau_n)\right),\quad \tau_i\in(0,h),i=1,\ldots,n.$$

The key observation which was the main motivation to develop rigorous numerical method based on this formula is that the coefficient $c_q^{q,p}=\binom{p+q}{q}^{-1}$ can be very small for $p=q$. Thus, this formula can have significantly smaller remainder than the Lagrange remainder used in the Taylor series method.

Now we would like to apply (\ref{eq:ho-formula}) to the flows $\varphi$ and $\psi:=D_x\varphi$. Let $[x_k]$ be a set of initial conditions and assume that from the predictor step we have computed $([x_{k+1}^0],[r^0],[V^0],[R^0])$ satisfying (\ref{eq:predictor-property1}--\ref{eq:predictor-property2}).

Let us fix positive integers $p,q$ such that $m=p+q$, $x_k\in [x_k]$ and put $x_{k+1}=\varphi(h,x_k)$. The formula (\ref{eq:ho-formula}) applied to this case reads
\begin{equation*}
 \sum_{i=0}^qc_i^{q,p}(-h_k)^i \tc \varphi i(0,x_{k+1}) =
  \sum_{i=0}^pc_i^{p,q}h_k^i \tc \varphi i(0,x_k) + \varepsilon,
\end{equation*}
where $\varepsilon \in (-1)^qc_q^{q,p}[r^0]$. Identifying vectors $x_k$, $x_{k+1}$ with unique solutions $x_k(\cdot)$, $x_{k+1}(\cdot)$ to the ODE passing through them at time zero, we obtain the equivalent but shorter form
\begin{equation}\label{eq:psiho-formula}
\Psi_{q,p}(-h_k,x_{k+1},0) = \Psi_{p,q}(h_k,x_k,0) + \varepsilon.
\end{equation}

Take the midpoints $\widehat x_{k+1}^0\in[x_{k+1}^0]$, $\widehat x_k\in[x_k]$. 
Since interval vectors are convex sets, and the local flow is a smooth function 
in both variables, we can apply the mean-value form to both sides of 
(\ref{eq:psiho-formula}) to obtain
\begin{equation*}
\Psi_{q,p}(-h_k,\widehat x_{k+1}^0,0) +  J_-(x_{k+1}-\widehat x_{k+1}^0) = 
    \Psi_{p,q}(h_k,\widehat x_k,0) + J_+(x_k-\widehat x_k) + \varepsilon
\end{equation*}
for some 
\begin{eqnarray*}
    J_- &\in& \left[D_x\Psi_{q,p}(-h_k,[x_{k+1}^0],0) \right]\quad\text{and}\\
    J_+ &\in& \left[D_x\Psi_{p,q}(h_k,[x_k],0)\right].
\end{eqnarray*}

We obtained a linear equation for $x_{k+1}$
\begin{equation}\label{eq:ho-implicit-equation}
    J_- (x_{k+1} - \widehat x_{k+1}^0) = J_+(x_k-\widehat 
x_k)+(\Psi_{p,q}(h_k,\widehat x_k,0) - \Psi_{q,p}(-h_k,\widehat 
x_{k+1}^0,0))+\varepsilon
\end{equation}
in which the matrices $J_\pm$ are unknown, but they can be rigorously bounded. 
Denoting
\begin{eqnarray*}
    [\delta] &=& \Psi_{p,q}(h_k,\widehat x_k,0) - \Psi_{q,p}(-h_k,\widehat 
x_{k+1}^0,0),\\
    \mbox{}[\varepsilon] &=& (-1)^qc^{q,p}_{q}[r^0],\\
    \mbox{}[J_-] &=& \left[D_x\Psi_{q,p}(-h_k,[x_{k+1}^0],0)\right],\\
    \mbox{}[J_+] &=& \left[D_x\Psi_{p,q}(h_k,[x_k],0)\right],\\
    \mbox{}[S] &=& \Id-\widehat J_-^{-1} [J_-],\\
    \mbox{}[r] &=& \widehat J_-^{-1}\left([\delta]+[\varepsilon]\right) + 
[S]([x_{k+1}^0]-\widehat x_{k+1}^0)
\end{eqnarray*}
and applying the interval Krawczyk operator 
\cite{Alefeld1994,Krawczyk1969,Neumaier1990} to the linear system 
(\ref{eq:ho-implicit-equation}), we obtain that for $x_k\in[x_k]$,
\begin{equation}\label{eq:ho-eval}
\varphi(h_k,x_k) = x_{k+1} \in  \widehat x_{k+1}^0 + \left(\widehat J_-^{-1} 
[J_+]\right)([x_k]-\widehat x_k) + [r]
\end{equation}
which is the main evaluation formula in the interval Hermite-Obreshkov method 
for IVPs presented in \cite{NedialkovJackson1998}. Note that this formula has 
exactly the same structure as (\ref{eq:predictor-eval}) used in the predictor 
step.

Each coefficients of $[S]$ is an interval containing 
zero and its diameter tends to zero with $h_k\to0$. The 
vector $[\delta]$ is almost a point vector, and 
$[\varepsilon]$ can be made as small as we need (manipulating the time step 
$h_k$). Therefore, the total error accumulated in $[r]$ is usually very thin in 
comparison to the size of the set $[x_k]$ we propagate. Thus, the main source 
of 
overestimation when evaluating (\ref{eq:ho-eval}) comes from the propagation of 
the product $\left(\widehat J_-^{-1} [J_+]\right)([x_k]-\widehat x_k)$. There is 
a wide 
literature on how to reduce this wrapping effect for such propagation (see 
\cite{MrozekZgliczynski2000} for a survey), and we will give some details 
concerning this issue in Section~\ref{sec:propagation}.

In what follows we argue that, with a little additional cost, we can compute a 
possibly tighter enclosure for the solutions to variational equation than the 
bound $[V^0]$ obtained from the predictor step. Let us fix $x_k\in[x_k]$ and 
let $V=\psi(h_k,x_k,\Id)$. Applying (\ref{eq:ho-formula}) to the solutions to 
the variational equation, we obtain
\begin{equation*}
 \sum_{i=0}^qc_i^{q,p}(-h_k)^i \tc \psi i(0,x_{k+1},V) =
  \sum_{i=0}^pc_i^{p,q}h_k^i \tc \psi i(0,x_k,\Id) + E,
\end{equation*}
where $E\in(-1)^qc_q^{q,p}[R^0]$. Since $\psi$ is linear in $V$, we obtain that 
the matrix $V=\psi(h_k,x_k,\Id)$ belongs to the solution set to the linear 
equation
\begin{equation*}
 [J_-]V = [J_+]+[E],
\end{equation*}
where
$[E]=(-1)^qc_q^{q,p}[R^0]$. Note that from the predictor step we already know 
that $V\in [V^0]$. Applying the interval Krawczyk operator 
\cite{Alefeld1994,Krawczyk1969,Neumaier1990} 
to this linear system we obtain
\begin{equation*}
 V\in \widehat J_-^{-1}([J_+]+[E]) + (\Id-\widehat J_-^{-1} [J_-])[V^0] =  
\widehat 
J_-^{-1}([J_+]+[E]) + [S][V^0].
\end{equation*}
Due to linearity of the variational equation, we can reuse the matrices 
$[J_-]$, $[J_+]$, $[S]$ and $\widehat J_-^{-1}$ computed in the corrector step 
for $\varphi$. 
Thus the additional cost is just a few matrix additions and multiplications.
Algorithm~\ref{alg:corrector} and Lemma~\ref{lem:alg2} summarize the above 
considerations.
\begin{algorithm}\label{alg:corrector}
 \SetKwInOut{Input}{Input}
 \SetKwInOut{Output}{Output}
 \Input{$p,q$ - positive integers\\
 $h_k$ - positive real number\\
 $[x_k]$ - interval vectors\\
 $([x_{k+1}^0],[r^0],[V^0],[R^0])$ - from the predictor step with\\
 $m=p+q$}
 \Output{$([x_{k+1}],[V])$}
 \textbf{Compute}:\\
  $[\delta]\gets \Psi_{p,q}(h_k,\widehat x_k,0) - \Psi_{q,p}(-h_k,\widehat 
x_{k+1}^0,0)$\;
  $[\varepsilon] \gets (-1)^qc^{q,p}_{q}[r^0]$\;
  $[J_-] \gets \left[D_x\Psi_{q,p}(-h_k,[x_{k+1}^0],0)\right]$\;
  $[J_+] \gets \left[D_x\Psi_{p,q}(h_k,[x_k],0)\right]$\;
  $[S] \gets \Id-\widehat J_-^{-1} [J_-]$\;
  $[r] \gets \widehat J_-^{-1}\left([\delta]+[\varepsilon]\right) + 
[S]([x_{k+1}^0]-\widehat x_{k+1}^0)$\;
  $[R]\gets \widehat J_-^{-1}\left([J_+]+(-1)^qc_q^{q,p}[R^0]\right)$\;
  $[x_{k+1}]\gets \left(\widehat x_{k+1}^0 + \left(\widehat J_-^{-1} 
[J_+]\right)([x_k]-\widehat x_k) + [r]\right)\cap[x_{k+1}^0]$\;
  $[V]\gets \left([R]+[S][V^0]\right)\cap[V^0]$\;
  \Return{$([x_{k+1}],[V])$}\;
\caption{Corrector.}
\end{algorithm}
\begin{lemma}\label{lem:alg2}
Assume that $h_k$, $[x_k]$, $[\widetilde y]$, $[\widetilde V]$ are such that 
(\ref{eq:c0enc}) and (\ref{eq:c1enc}) hold and that the quadruple 
$([x_{k+1}^0],[r^0],[V^0],[R^0])$ is returned by the predictor step 
(Algorithm~\ref{alg:predictor}). Then the quantities $([x_{k+1}],[V])$ computed 
by Algorithm~\ref{alg:corrector} satisfy
\begin{eqnarray}
 \varphi(h_k,[x_k]) & \subset & [x_{k+1}]\label{eq:corrector-property1} 
\quad\textrm{ and }\\
 \psi(h_k,[x_k],\Id) & \subset & [V] \label{eq:corrector-property2}.
\end{eqnarray}
\end{lemma}

We would like to emphasize that by its construction the proposed algorithm always returns tighter bounds than the \lo algorithm because the result obtained from the corrector step is intersected with the bound obtained from the predictor step.

\subsection{Propagation of product of interval objects.}\label{sec:propagation}
It is well known that evaluation of expressions in interval arithmetic can 
produce large overestimation due to dependency of variables and the wrapping 
effect \cite{Alefeld1994,Lohner1992,Moore1966,Neumaier1990}. To reduce this 
undesirable drawback we follow the ideas from 
\cite{Lohner1992,MrozekZgliczynski2000,NedialkovJackson1998,Zgliczynski2002}, 
and we represent subsets of $\mathbb R^n$ and $\mathbb R^{n^2}$ in the forms 
(doubletons in \cite{MrozekZgliczynski2000} terminology)
\begin{eqnarray}
\mbox{}[x_k] = x_k + C_k[r_k] + 
B_k[s_k]\label{eq:xRepresentation}\quad\text{and}\\
\mbox{}[V_k] =  V_k + A_k[R_k] + Q_k[S_k]\label{eq:VRepresentation}.
\end{eqnarray}
The initial conditions $([x_0],[V_0])$ of (\ref{eq:vivp}) are assumed to be already in the form (\ref{eq:xRepresentation}--\ref{eq:VRepresentation}). The parallelepipeds $x_k + C_k[r_k]$ and $V_k + A_k[R_k]$ are used to store the main part of the sets $[x_k]$ and $[V_k]$, respectively. The terms $B_k[s_k]$ and $Q_k[S_k]$ are used to collect all usually thin quantities that appear during the computation.

According to (\ref{eq:V-eval}), the bound for $\psi(t_k+h_k,[x_0],[V_0])$ can be computed as
\begin{equation*}
[V_{k+1}] = [V][V_k],
\end{equation*}
where $[V]$ satisfies (\ref{eq:corrector-property2}). Substituting the 
representation (\ref{eq:VRepresentation}) we obtain
\begin{equation*}
 [V_{k+1}] \subset [V]\left(V_k + A_k[R_k] + Q_k[S_k]\right)\cap \left(V_{k+1} + A_{k+1}[R_{k+1}] + Q_{k+1}[S_{k+1}]\right),
\end{equation*}
where the new representation is computed as follows
\begin{eqnarray*}
 \mbox{}[\Delta A]&=&\left([V]-\widehat V\right)\left(V_k + A_k[R_k]\right),\\
  V_{k+1}&=& \widehat VV_k,\\
  A_{k+1}&=& \widehat V A_k,\\
 \mbox{}[S_{k+1}]&=& \left(Q_{k+1}^{-1}[V][Q_k]\right)[S_k] + 
Q_{k+1}^{-1}[\Delta A]\quad\text{and}\\
 \mbox{}[R_{k+1}]&=&[R_k].
\end{eqnarray*}
In principle, the matrix $Q_{k+1}$ can be chosen as any invertible matrix. The 
numerical experiments 
\cite{Lohner1992,MrozekZgliczynski2000,NedialkovJackson1998,Zgliczynski2002} 
show that one of the most efficient strategies in reducing the wrapping effect 
is to compute $Q_{k+1}$ as an orthogonal matrix from the $QR$ decomposition of 
the point matrix $\widehat V Q_k$. Note, that even if the matrix $Q_{k+1}$ is a 
point matrix, the inverse $Q_{k+1}^{-1}$ must be computed rigorously in interval 
arithmetic.

Similar strategy is used for propagation of products in
\begin{eqnarray*}
 [x_{k+1}] & \subset & \widehat x_{k+1}^0 + \left(\widehat J_-^{-1} 
[J_+]\right)([x_k]-\widehat x_k) + [r]\\
  & = & \widehat x_{k+1}^0 + \left(\widehat J_-^{-1} [J_+]\right)(C_k[r_k] + 
B_k[s_k]) + [r]
\end{eqnarray*}
--- see 
\cite{Lohner1992,MrozekZgliczynski2000,NedialkovJackson1998,Zgliczynski2002} for 
details. 
\section{Complexity.}\label{sec:complexity}
In this section, we explain why the \ho algorithm may perform better than the \lo algorithm, even if it has higher computational complexity. A large numerical and theoretical study were performed to compare the Interval Hermite-Obreshkov method (IHO) with the Interval Taylor Series Method (ITS) \cite{NedialkovJacksonCorliss1999}. It has been shown that, with the same step size and order, the IHO method is more stable and produces smaller enclosures than the ITS method on constant coefficient problems. Furthermore, the IHO method allows the use of a much larger stepsize than the ITS method, thus saving computation time during the whole integration. However, comparing these two methods in the nonlinear case is not as simple as in the constant coefficient case. Our goal is to predict the benefits of performing additional calculations required by the IHO method applied to (\ref{eq:vivp}). 

\subsection{Cost of \lo and \ho methods per step.}
We assume that both predictor (Algorithm~\ref{alg:predictor}) and corrector 
(Algorithm~\ref{alg:corrector}) have the same order. That is, if the order of 
the predictor is $m$, we consider the corrector step with $p$ and $q$ such that 
$m=p+q$. In what follows we list the most time-consuming items of the predictor 
and corrector, which are the core of their computational complexity. 

In the analysis give below, we count the number of operations which are really 
executed by the implementation, rather than the possible theoretical and 
asymptotic complexity. Therefore, we assume that the product of two square 
interval matrices is computed by the naive algorithm (three nested loops or 
equivalent), which executes exactly $n^3$ interval multiplications. 

We would like to emphasize, that the rigorous integration of a differential 
equation is a very difficult task even in quite low dimensions. Thus, 
dimensions used in practice are usually less than $20$. 
Computer-assisted proofs for $100$-dimensional systems are actually the state 
of the art --- see for instance \cite{KapelaSimo2007}. Therefore, the use 
of asymptotically fast algorithms for matrix 
multiplications, such as the Strassen algorithm \cite{Strassen} or the 
Coppersmith-Winograd 
\cite{COPPERSMITH1990251} algorithm, does not make any sense.

Let us denote by $c_f$ the cost of evaluating the vector field (\ref{eq:vivp}). 
For the \lo step we need the following operations (predictor step and 
propagation of doubleton representations)
\begin{itemize}
 \item simultaneous computation of $\tc \varphi i (0,[x_k])$ and $\tc \psi i 
(0,[x_k],\Id)$ up to order $m$. This is performed by means of automatic 
differentiation techniques, and it takes $c_f(2n+1)(m+1)(m+2)/2$ 
multiplications --- see \cite{RallCorliss1996},
 \item simultaneous computation of $\tc \varphi i (0,[\widetilde y])$ and $\tc 
\psi 
i (0,[\widetilde y],\Id)$ up to order $m+1$. This is performed by means of the 
automatic differentiation techniques and it takes $ 
c_f(2n+1)(m+2)(m+3)/2$ multiplications --- see \cite{RallCorliss1996},
 
 \item $13$ matrix by matrix multiplications, $2$ point matrix inversions and 
$2$ point matrix QR decompositions. Approximate QR decomposition of a point 
matrix is much cheaper than the product of interval matrices and we may assume 
that it takes $O(n^3)$ with a constant less than one (in terms of interval 
multiplications).  The inversion of a point matrix which is very close to 
orthogonal is performed by means of the interval Krawczyk operator 
\cite{Krawczyk1969} and takes $n^3$ (one multiplication) because an approximate 
result is already known (transposition of an approximate orthogonal matrix). 
Thus, the total cost of all matrix operations listed above is at most $17n^3$.
\end{itemize}
We did not list cheaper operations like additions, intersections of interval 
objects, matrix by vector products. All polynomial evaluations perform in 
total $O(n^2m)$ interval multiplications, and they add significant cost to 
linear systems ($c_f=0$) and to nonlinear systems but with very small number 
of nonlinear terms ($c_f \ll n$). Thus, we skip them.

To sum up, the total costs of the \lo step is 
\begin{equation*}
C_{\textrm{LO}}(n,m)\simeq c_f(2n+1)(m+2)^2 +17n^3.
\end{equation*}

In the \ho method, we can reuse the Taylor coefficients of $\varphi$ and $\psi$ 
computed in the predictor step, which are needed for computing the $[J_+]$ 
matrix and $\Psi_{p,q}(h_k,\widehat x_k,0)$. Thus, the additional cost is
\begin{itemize}
 \item computation of $\tc \psi i (0,[x_{k+1}^0],\Id)$ up to order $q$. This is 
performed by means of automatic differentiation techniques, and it takes 
$c_f(2n+1)(q+1)(q+2)/2$ operations --- see \cite{RallCorliss1996},
 \item computation of $\Psi_{q,p}(-h_k,\widehat x_{k+1}^0,0)$ takes 
$c_f(q+1)(q+2)/2$,
 \item rigorous inversion of the point matrix $\widehat J_-$ takes at most 
$2n^3$ (one nonrigorous inverse and one interval matrix multiplication in 
the Krawczyk method) and 
\item $4$ interval matrix multiplications require in total $4n^3$ operations.
\end{itemize}

The total additional cost of the \ho step is at most
\begin{equation*}
c_f(n+1)(q+1)(q+2) + 6n^3.
\end{equation*}

\begin{figure}[htbp]
 \centerline{\includegraphics[width=\textwidth]{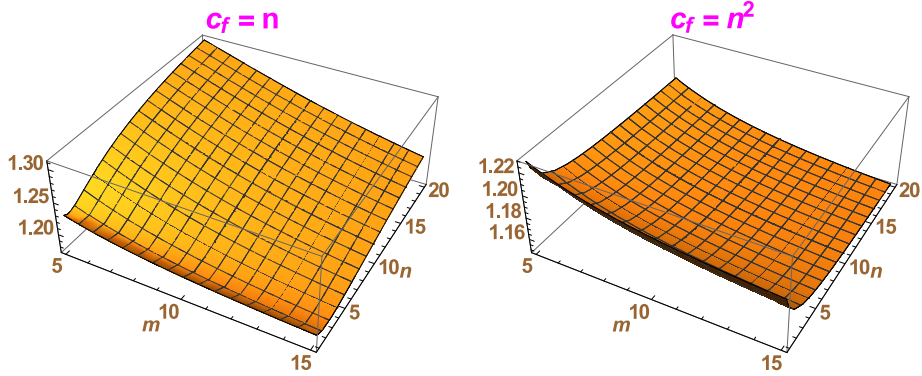}}
 \caption{Plot of $C_{\textrm{HO}}(n,m)/C_{\textrm{LO}}(n,m)$ for $c_f=n$ and $c_f=n^2$, 
respectively.\label{fig:cost}}
\end{figure}

Assume that $m$ is an even number and take $q=p=\frac{m}{2}$. Then the above 
additional cost of the \ho method is approximately
\begin{equation*}
\frac{1}{4}c_f(n+1)(m+2)(m+4) + 6n^3.
\end{equation*}
Hence, total computational complexity of the \ho method is
\begin{equation*}
C_{\textrm{HO}}(n,m)\simeq C_{\textrm{LO}}(n,m)+\frac{1}{4}c_f(n+1)(m+2)(m+4) + 
6n^3.
\end{equation*}
In general, the complexity depends on the cost of the vector field evaluation 
$c_f$ which can be arbitrary. In Fig.~\ref{fig:cost} we plot the graph of 
$C_{\textrm{HO}}/C_{\textrm{LO}}$ for two cases. The case $c_f=n$ means that the number of 
nonlinear terms in the vector field is equal to the dimension of the problem. We 
observe that, if order $m$ of the method is much smaller than the dimension $n$, 
then the complexity is dominated by the matrix operations and we have
\begin{equation*}
 \lim_{n\to\infty} C_{\textrm{HO}}(n,m)/C_{\textrm{LO}}(n,m) = \frac{23}{17}\approx 1.35294.
\end{equation*}
for all fixed values of $m$. We observe, however, that for reasonable 
 dimensions and orders, this factor is much smaller 
than the limit value.

A model example for the $c_f=n^2$ case is a second order polynomial vector field 
with nonzero coefficients in the quadratic terms. In this case we have
\begin{equation*}
 \lim_{n\to\infty} C_{\textrm{HO}}(n,m)/C_{\textrm{LO}}(n,m) = \frac{9 m^2+38 m+132}{8 
m^2+3m+100}.
\end{equation*}

The above analysis shows that the additional cost of the \ho method in a typical 
nonlinear case approaches $1/8$. In the next section, we argue that this extra 
cost of the \ho method is compensated by the larger time steps this method can 
perform without losing the accuracy.

\subsection{Maximal allowed time step for a fixed error tolerance.}
To obtain insights into the compared methods, we ask the following question: 
given an acceptable tolerance $\varepsilon$ per step, what is the maximal time 
step $h$ of both methods that guarantees achieving this constraint. For the \lo 
method, we have to solve the following inequality
\begin{equation*}
\left \|h^{m+1} {\tc \varphi {m+1} (0,[\widetilde y])} \right\| \leq 
\varepsilon.
\end{equation*}
In general, it is very difficult to answer this question because $[\widetilde 
y]=[\widetilde y(h)]$ depends on $h$. If $[x_k]$ is a point and $\varepsilon$ is 
very small, we can assume that the vector field is almost constant near $[x_k]$ 
and thus $\tc \varphi {m+1} (0,[\widetilde y])\approx \tc \varphi {m+1} 
(0,[x_k])$. Since $[x_k]\subset [\widetilde y]$, we always have $\|\tc \varphi 
{m+1} (0,[x_k])\|\leq \|\tc \varphi {m+1} (0,[\widetilde y])\|$. With this 
simplification, we obtain an upper bound for the time step 

\begin{equation*}
h_{\textrm{LO}}:= h = \sqrt[m+1]{\frac{\varepsilon}{\left\|\tc \varphi {m+1} 
(0,[x_k]) \right\|}}.
\end{equation*}

\begin{figure}[htbp]
 \centerline{\includegraphics[width=.7\textwidth]{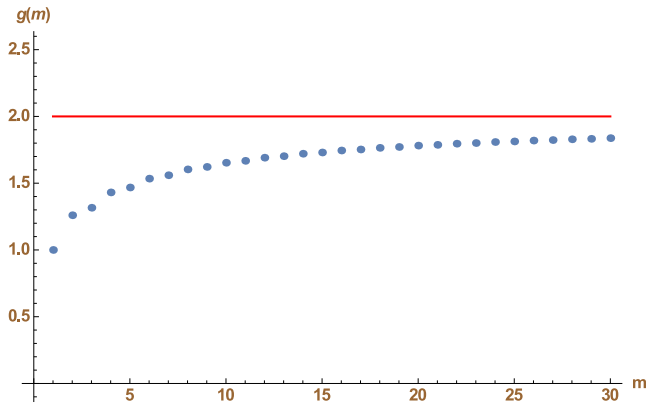}}
 \caption{Plot of the theoretical maximal factor of maximal time step in the \ho 
and \lo methods for a fixed tolerance --- see 
\ref{eq:functiongm}.\label{fig:stepsize}}
\end{figure}

For the \ho method we obtain the following upper bound for the time step 
\begin{equation*}
h_{\textrm{HO}}:= h = 
\sqrt[m+1]{\binom{m}{{\lceil\frac{m}{2}}\rceil}\frac{\varepsilon}{\left\|\tc 
\varphi {m+1} (0,[x_k]) \right\|}},
\end{equation*}
where by $\lceil m/2\rceil$ we denote the smallest integer not smaller than 
$m/2$. Denote
\begin{equation}\label{eq:functiongm}
g(m):=h_{\textrm{HO}}/h_{\textrm{LO}} = 
\sqrt[m+1]{\binom{m}{\lceil\frac{m}{2}\rceil}}.
\end{equation}
It is easy to show that 
\begin{equation*}
\lim_{m\to\infty}g(m) = 2.
\end{equation*}
In Fig.~\ref{fig:stepsize}, we observe that the values of $g(m)$ rapidly grow 
for small values of $m$. This is important from practical point of view --- even 
for small order $m=6$ the \ho method allows up to 53\%  larger time steps than 
the \lo method. For $m=16$ this is 74\%. For larger values of the tolerance 
$\varepsilon$, the computed enclosure $[\widetilde y]$ for $h=h_{\textrm{LO}}$ 
is usually significantly smaller than that computed for $h=h_{\textrm{HO}}$, 
which affects the norm $\left\|\tc \varphi {m+1} (0,[\widetilde y])\right\|$. 
Therefore, the value $g(m)$ is a theoretical upper bound for the possible growth 
ratio of the time step in the \ho method achievable when $\varepsilon\to 0$. 

\section{Benchmarks.}\label{sec:benchmarks}
In this section, we present the results of a comparison of the \lo algorithm and the \ho algorithm. The structure of the tests is as follows. For a given ODE 
\begin{itemize}
 \item we take an initial condition $u$ which is an approximate periodic orbit for the system;
 \item we integrate the variational equations along this periodic orbit using 
the \lo and \ho algorithms with the same algorithm for rough enclosure (HOE), 
the same order  $m=p+q$ of the methods and a constant time step $h$;
 \item we compare the logarithm of the maximal diameter of the interval matrix $[V_k]$ (diameter of the widest component) computed by means of the two algorithms; and
 \item we repeat the above two steps six times: for two different orders of the numerical methods each for three different time steps.
\end{itemize}
Fixing the time steps allows us to compare the size of the enclosures returned by the two algorithms over the same time step. This will allow us to conclude that the \ho algorithm can take larger time steps than the \lo algorithm without significant lost of accuracy. The comparison of the two algorithms with variable time steps will be given in Section~\ref{sec:applications}.

The above test is performed for four ODEs: the Lorenz \cite{Lorenz1963} system, 
the H\'enon-Heiles system \cite{HenonHeiles1964}, the Planar Circular 
Restricted Three Body Problem (PCR3BP), and a $10$-dimensional moderately stiff 
ODE. Below we give initial conditions and discuss obtained results. 

The Lorenz system \cite{Lorenz1963} for ``classical'' parameters is given by
\begin{equation}\label{eq:lorenz}
\left\{
    \begin{array}{rcl}
        \dot x &=& 10(y-x),\\
        \dot y &=& x(28-z)-y,\\
        \dot z &=& xy-\frac{8}{3}z.
    \end{array}
\right.
\end{equation}
 
The H\'enon-Heiles system \cite{HenonHeiles1964} is a hamiltonian ODE given by 
\begin{eqnarray}\label{eq:hh}
\left\{
\begin{array}{lcl}
\ddot x &=& -x(1+2y),\\
\ddot y &=& x^2 - y(1+y).
\end{array}\right.
\end{eqnarray}

The PCR3BP is a mathematical model that describes motion of a small body with negligible mass in the gravitational influence of two big bodies. The motion is restricted to the plane, and the two main primaries rotate around their common mass centre. The equations for motion of the small body is then given by
\begin{equation}\label{eq:pcr3bp}
\begin{cases}
\ddot{x} - 2\dot y = D_x\Omega(x,y),\\
\ddot{y} + 2\dot x = D_y\Omega(x,y),
\end{cases}
\end{equation}
where
\begin{equation*}
\Omega(x,y) = \frac{1}{2}(x^2+y^2) + \frac{1-\mu}{\sqrt{(x+\mu)^2+y^2}} + \frac{\mu}{(x-1+\mu)^2+y^2}.
\end{equation*}
The parameter $\mu$ stands for the relative mass of the two main bodies. For our 
tests we fixed $\mu=0.0009537$, which corresponds to the Sun-Jupiter system.

The last ODE is the Galerkin projection of the following infinite dimensional ODE
\begin{equation}\label{eq:ks}
\dot a_k = k^2(1-\nu k^2)a_k - k \sum_{n=1}^{k-1}a_na_{k-n} + 2k\sum_{n=1}^\infty a_n a_{n+k}
\end{equation}
onto $(a_1,\ldots,a_{10})$ variables. The above system describes solutions to the one-dimensional Kuramoto-Sivashinsky PDE \cite{KuramotoTsuzuki1976,Sivashinsky1977} under periodic and odd boundary conditions, see \cite{Zgliczynski2004,ZgliczynskiMischaikow2001} for derivation.

\begin{figure}[htbp]
 \centerline{\includegraphics[width=\textwidth]{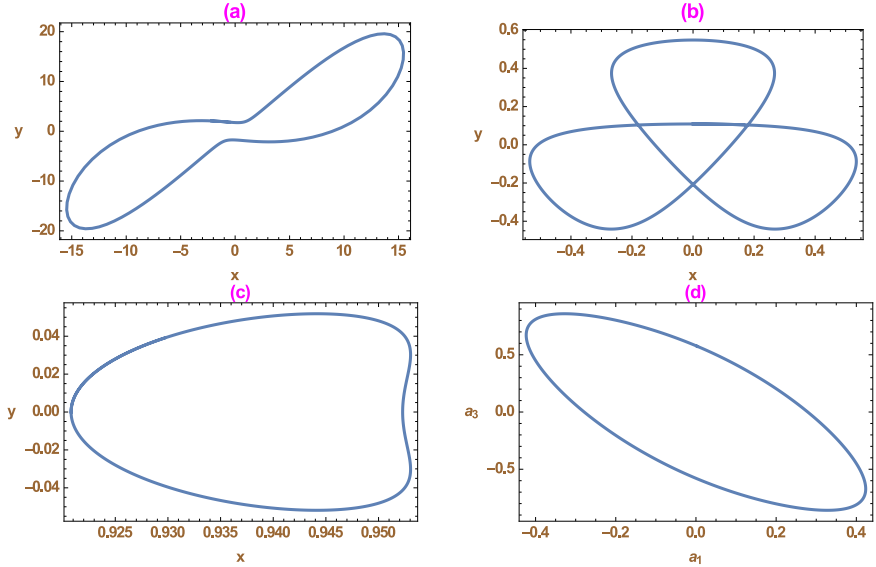}}
 \caption{Approximate periodic orbits for (a) the Lorenz system 
(\ref{eq:lorenz}), (b) H\'enon-Heiles hamiltonian (\ref{eq:hh}), (c) 
the PCR3BP (\ref{eq:pcr3bp}) and (d) a the $10$-dimensional Galerkin projection 
of the Kuramoto-Sivashinsky equation (\ref{eq:ks}), 
respectively.\label{fig:orbits}}
\end{figure}

We have chosen initial conditions that are close to periodic orbits of these systems (see Fig.~\ref{fig:orbits})
\begin{eqnarray*}
u_{\text{Lorenz}}&=&(-2.1473681756955529387,2.078047612582596404,27),\\
u_{\text{H\'enon-Heiles}}&=&(0.0,0.10903,,0.5677233993382853,0.0),\\
u_{\text{PCR3BP}}&=&(0.92080349132074,0.0,0.0,0.1044476727069111)\quad\text{
and}\\
u_{\text{KS}}&=&\left[\begin{matrix}    
    0.2012106\\
    1.2899797585174486\\
    0.2012106\\
    -0.37786628185377774\\
    -0.042309451521292417\\
    0.043161614695331821\\
    0.0069402112803455653\\
    -0.0041564870501656455\\
    -0.00079448972725675504\\
    0.00033160609117820303\end{matrix}\right].
\end{eqnarray*}
For the two Hamiltonian systems, the coordinates are given in the order $(x,y,\dot x, \dot y)$. The orbit $u_{\text{PCR3BP}}$ is the well known $L_1$ Lyapunov orbit for the Sun-Jupiter-Oterma system. In \cite{Zgliczynski2004}, a computer assisted proof of the existence of a periodic solution for the full infinite dimensional system (\ref{eq:ks}) is given. The projection of this periodic orbit onto the first $10$ coordinates is very close to the point $u_{\text{KS}}$. In fact, due to very strong dissipation, the variables with high indexes have very small impact on the dynamics of (\ref{eq:ks}). The system becomes very stiff even for relative small dimension of the Galerkin projection. 

\begin{figure}[htbp]
 \centerline{\includegraphics[width=\textwidth]{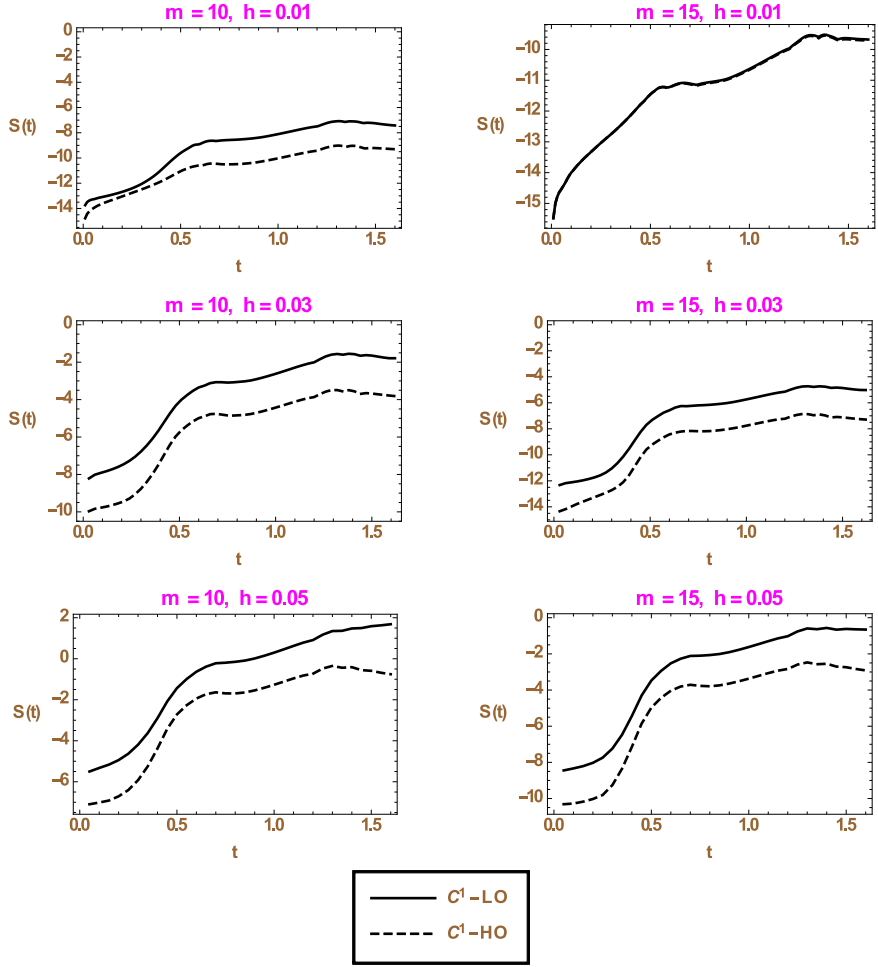}}
 \caption{The results of the tests for the Lorenz system. We plot 
$S(t)=\log_{10} \mathrm{diam}([V(t)])$ along trajectory of the point 
$u_{\text{Lorenz}}$ integrated with order $m$ and with fixed time step 
$h$.\label{fig:lorenz}}
\end{figure}

\begin{figure}[htbp]
 \centerline{\includegraphics[width=\textwidth]{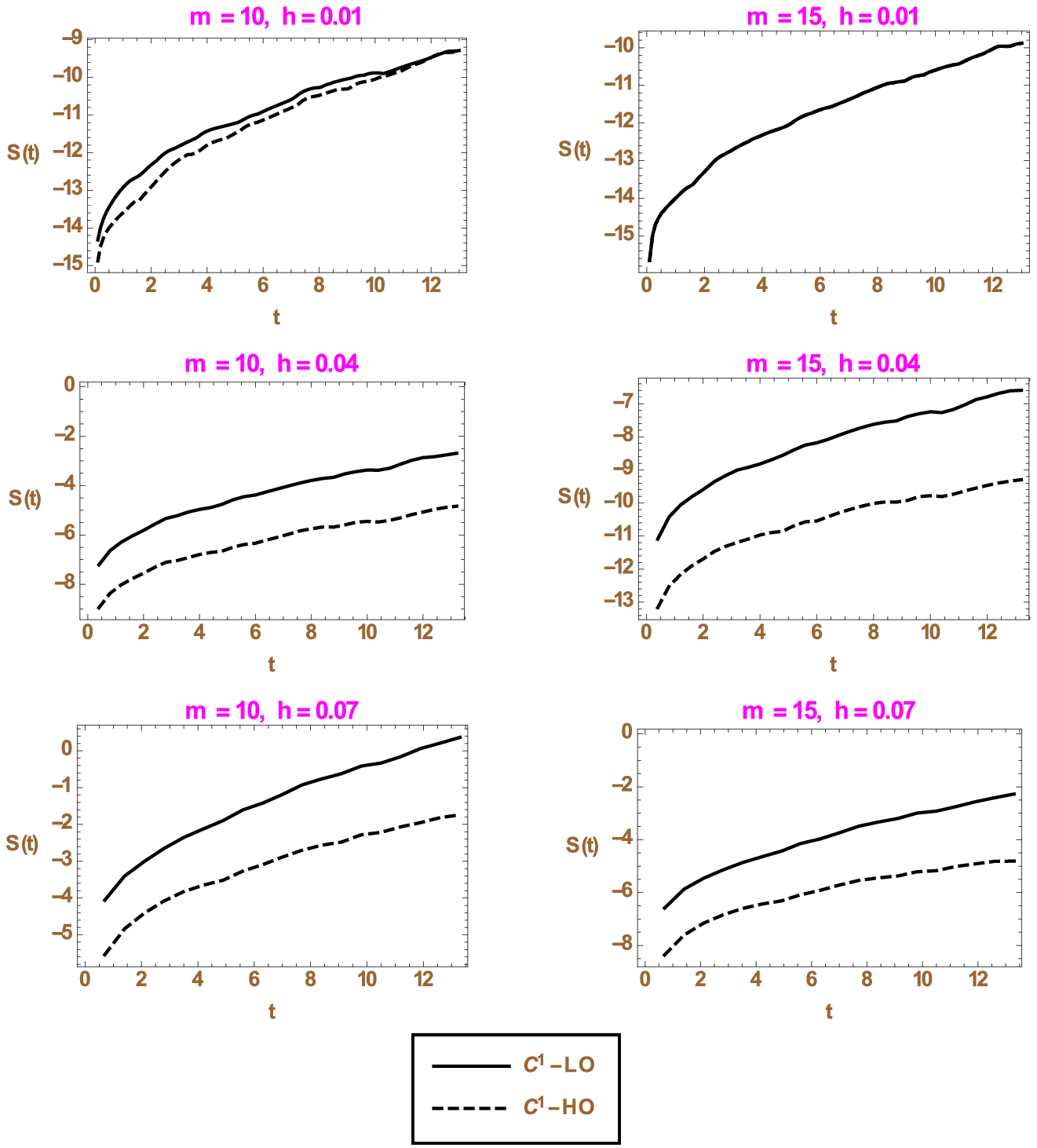}}
 \caption{The results of the tests for the H\'enon-Heiles system. We plot 
$S(t)=\log_{10} \mathrm{diam}([V(t)])$ along trajectory of the point  
$u_{\text{H\'enon-Heiles}}$ integrated with order $m$ and with fixed time step 
$h$.\label{fig:hh}}
\end{figure}

\begin{figure}[htbp]
 \centerline{\includegraphics[width=\textwidth]{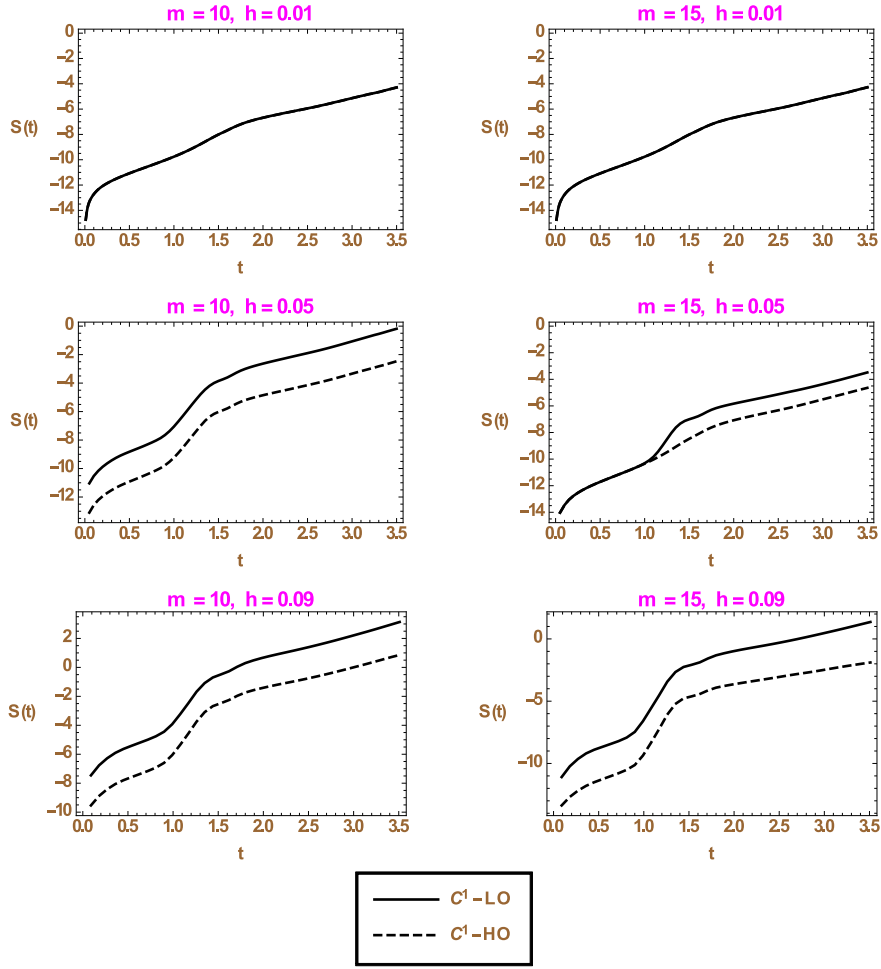}}
 \caption{The results of the tests for the PCR3BP. We plot 
$S(t)=\log_{10} \mathrm{diam}([V(t)])$ along 
trajectory of the point $u_{\text{PCR3BP}}$ integrated with order $m$ and with 
fixed time step 
$h$.\label{fig:pcr3bp}}
\end{figure}

\begin{figure}[htbp]
 \centerline{\includegraphics[width=\textwidth]{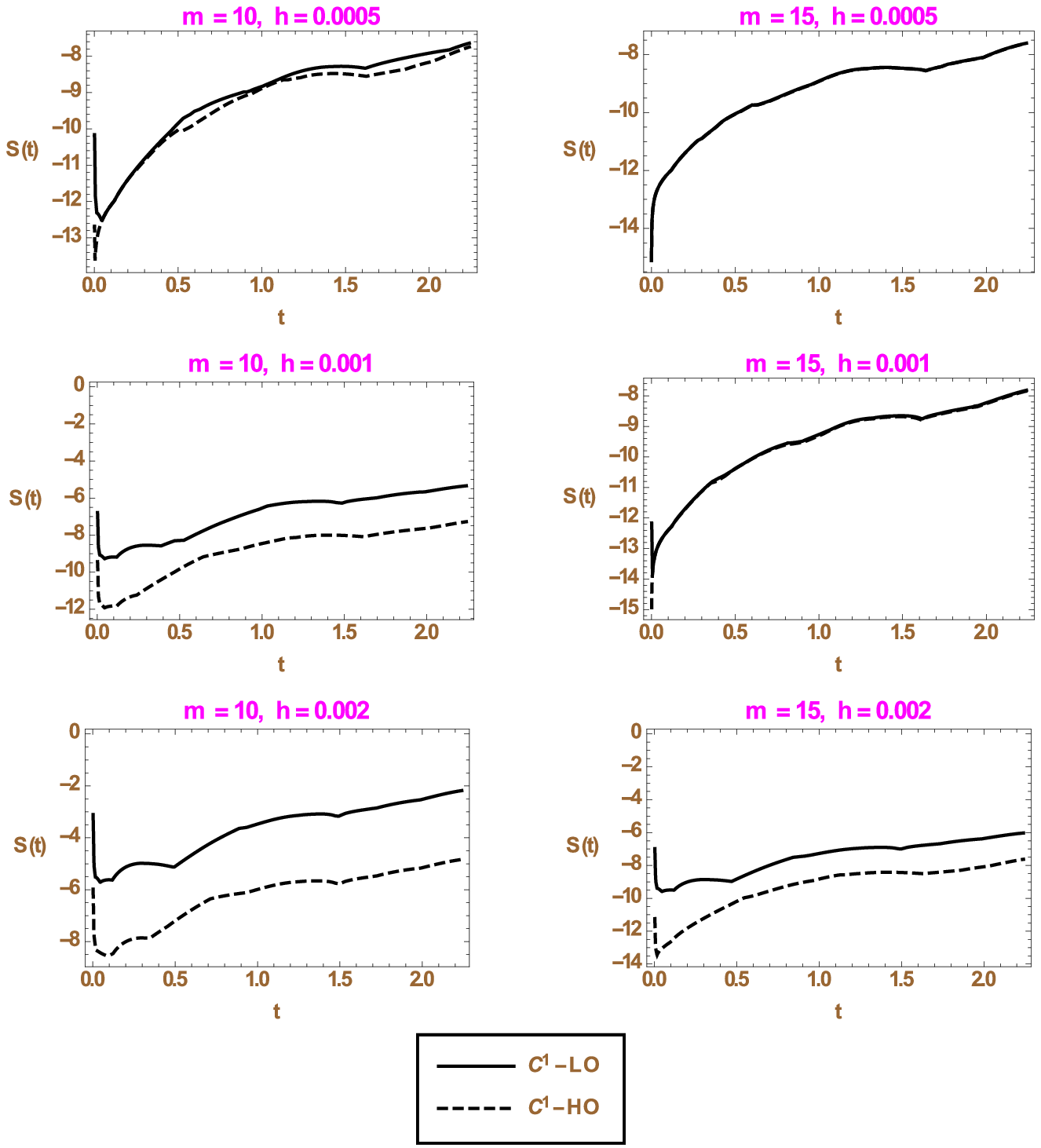}}
 \caption{The results of the tests for the 10-dimensional Galerkin projection 
of 
the Kuramoto-Sivashinsky equation. We plot 
$S(t)=\log_{10} \mathrm{diam}([V(t)])$ along  
trajectory of the point $u_{\text{KS}}$ integrated with order $m$ and with fixed 
time step 
$h$.\label{fig:ks}}
\end{figure}

In Figs.~\ref{fig:lorenz}--\ref{fig:ks} we present results of our numerical 
experiments. On these figures we show plot of $S(t)=\log_{10} 
\mathrm{diam}([V(t)])$ along an approximate periodic trajectory, where 
$\mathrm{diam}([V(t)])$ is the largest width of coefficient in the 
interval matrix $[V(t)]$. We can see that in each case the \ho method does not 
return worse results than the \lo algorithm. This is due to its construction, 
because the bounds computed in the corrector step are intersected with the 
estimates from the predictor step, which is used in the \lo algorithm. Indeed, 
in 
the Algorithm~\ref{alg:corrector} we have
$$[V]\gets \left([R]+[S][V^0]\right)\cap[V^0].$$
Looking at the columns of Figs.~\ref{fig:lorenz}, \ref{fig:hh} and \ref{fig:pcr3bp}, we observe that in each case the advantage of the \ho method increases when the time step is enlarged, and the obtained bounds can be orders of magnitude tighter. This is due to the fact that the \ho method has $c^{q,p}_{q}$ times tighter truncation error than the Taylor method. To give some numbers, let us take $m=20$ which is a typical order used in computations. Then $p=q=10$ and $c^{q,p}_{q}=c_{10}^{10,10}\approx 5.4\cdot 10^{-6}$.

Increasing the order of the method makes the truncation error of the \lo method 
smaller when the time step is fixed. Therefore, in the  right columns of each 
figure we observe that the \lo method performs much better than in the left 
column. The \ho method, however, still returns tighter enclosures and is capable 
to take even larger time steps without significant lost of accuracy. This is 
especially important for stiff problems, where the time steps used by a 
nonstiff solver cannot be large and thus integration over large time interval 
is very expensive. We would like to emphasize, that the maximal possible time 
step that a rigorous ODE solver can take is limited mainly by the 
possibility of finding a rough enclosure over the time step. To the best of our 
knowledge, the HOE algorithm \cite{NedialkovJacksonPryce2001}, which is 
nonstiff, is one of the most efficient. Therefore construction of a 
general rigorous stiff ODE solver without extra knowledge of the system is a 
challenge. Remarkable exceptions are solvers for infinite-dimensional 
strongly dissipative systems \cite{Cyranka2013,Zgliczynski2004}, where the 
structure of the system is used to construct a dedicated so-called dissipative 
enclosure.

In  Fig.~\ref{fig:ks} we can see that the \ho method can perform much larger 
time steps keeping very good accuracy of computed bounds.

The bounds obtained by the \ho method are tighter than those returned by the 
\lo algorithm, but as we observed in Section~\ref{sec:complexity}, the \ho 
method is computationally more expensive. In the next section, we argue that 
this extra cost per step is compensated by the larger time steps we can take.
\section{Applications.}\label{sec:applications}
In this section, we present an application of the proposed algorithm to a~
computer-assisted proof of a new result concerning the R\"ossler system 
\cite{Rossler1976}. We focus on the comparison of the time of computation needed 
to prove this result, when the \lo algorithm and the \ho algorithm are used to 
integrate variational equations, which are necessary to prove this theorem.

The classical R\"ossler system \cite{Rossler1976} is given by (\ref{eq:rosslerSystem}). When the two parameters $a,b$ vary, the system exhibits wide spectrum of bifurcations. This system admits period doubling bifurcations \cite{WilczakZgliczynski2009focm}, which lead to chaotic dynamics \cite{Zgliczynski1997}. In \cite{Pilarczyk2003}, the existence of two periodic orbits was proved by means of the Conley index theory. All the above results about the system (\ref{eq:rosslerSystem}) are computer assisted and use rigorous ODE solvers.

\begin{figure}[htbp]
\centerline{
  \includegraphics[width=.6\textwidth]{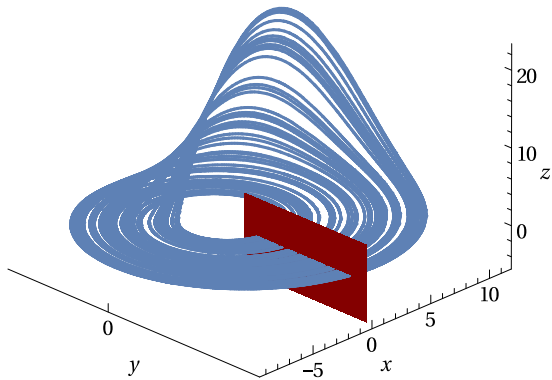}
}
\caption{Typical chaotic trajectory of the system (\ref{eq:rosslerSystem}) and a slice of the Poincar\'e section $\Pi$.\label{fig:PoincareSection}}
\end{figure}
Let $\Pi=\{(x,y,z)\in\mathbb R^3: x=0 \text{ and } \dot x>0\}$ be a Poincar\'e
section (see Fig.~\ref{fig:PoincareSection}) and let $P:\Pi\to\Pi$ be 
the Poincar\'e map. Since the $x$ coordinate is equal to zero on $\Pi$, we use 
only two coordinates $(y,z)$ to describe points on $\Pi$.

\begin{thm}\label{thm:rosslerMainResult}
 Let $l_B=-10.7$, $r_B = -2.3$, $l_M=-8.4$, $r_M = -7.6$, $l_N=-5.7$, $r_N=-4.6$, $Z=[0.028,0.034]$ and let
 \begin{eqnarray*}
  B &=& [l_B,r_B]\times Z,\\
  M &=& [l_M,r_M]\times Z \textrm{ and }\\
  N &=& [l_N,r_N]\times Z.
 \end{eqnarray*}
For the classical parameter values $a=0.2$, $b=5.7$ the following statements hold.
\begin{itemize}
\item The system (\ref{eq:rosslerSystem}) admits an attractor. The set $B$ is a trapping region for the Poincar\'e map, i.e. $P$ is well defined on $B$ and $P(B)\subset B$. In particular, there exists a maximal invariant set $\mathcal A = \bigcup_{n>0}P^n(B)$ for the map $P$ that is compact and connected.

\item The maximal invariant set for $P^2$ in $N\cup M$, denoted by $\mathcal H 
=\mathrm{inv}(P^2,N\cup M)\subset \mathcal A$, is uniformly hyperbolic; in 
particular it is robust under perturbations of the system. The dynamics of 
$P^2$ on $\mathcal H$ is chaotic in the sense that $P^2|_{\mathcal H}$ is 
conjugated to the Bernoulli shift on two symbols.
\end{itemize}
\end{thm}

\begin{figure}[htbp]
 \centerline{\includegraphics[width=.6\textwidth]{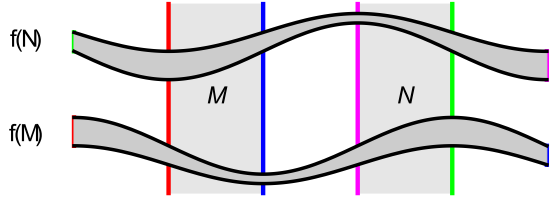}}
 \caption{Geometric conditions that guarantee the existence of chaotic dynamics 
--- see (\ref{eq:covrel}).\label{fig:covrel}}
\end{figure}

\textbf{Proof:}
The tools used in a computer-assisted proof of 
Theorem~\ref{thm:rosslerMainResult} are well known, and we summarize them here.

\begin{figure}[htbp]
 \centerline{\includegraphics[width=.6\textwidth]{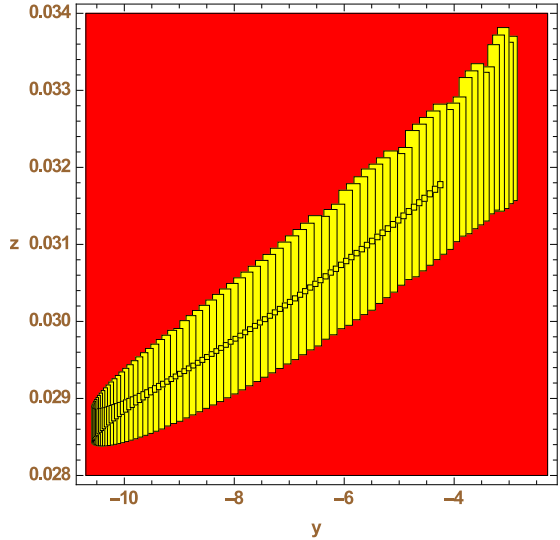}}
 \caption{The set $B$ (in red) and a rigorous enclosure for $P(B)$ (in 
yellow) obtained as the union of enclosures $\bigcup_{i=1}^{160} P(B_i)$ --- 
see (\ref{eq:trapping_region}).\label{fig:trapping_region}}
\end{figure}

\textbf{Trapping region}. Verification that $B$ is a trapping region for $P$ reduces to checking the inclusion
 \begin{equation*}
  P(B)\subset B.
 \end{equation*}
 We uniformly subdivided the set $B$ onto $N=160$ pieces of the form $B_i = 
[y_{i},y_{i+1}]\times Z$, $y_i=l_B + i \cdot (r_B-l_B)/N$. 
 Then we verified that 
 \begin{equation}\label{eq:trapping_region}
 \bigcup_{i=1}^N P(B_i)\subset B.
 \end{equation}
 We used a rigorous ODE solver of order $25$ from the CAPD 
library which implements the $\mathcal C^0$ Hermite-Obreshkov algorithm proposed 
in \cite{NedialkovJackson1998}. Rigorous enclosure for $P(B)$ returned by our 
routine is shown in Fig.~\ref{fig:trapping_region}.
 
\textbf{Chaos}. Semiconjugacy of $P^2|_{\mathcal H}$ to the Bernoulli shift is 
proved by means of the method of covering relations --- the same as in 
\cite{Zgliczynski1997} but applied to different sets. It is sufficient to check 
the following geometric conditions
  \begin{equation}\label{eq:covrel}
  \begin{array}{lclc}
  \pi_yP^2(y,z)&<& l_M \quad \text{for} (y,z)\in\{l_M\}\times Z,\\  
  \pi_yP^2(y,z)&>& r_N \quad \text{for} (y,z)\in\{r_M\}\times Z,\\  
  \pi_yP^2(y,z)&<& l_M \quad \text{for} (y,z)\in\{r_N\}\times Z & \textrm{ and 
} \\  
  \pi_yP^2(y,z)&>& r_N \quad \text{for} (y,z)\in\{l_N\}\times Z,  
  \end{array}
  \end{equation}
where $\pi_y$ denotes the canonical projection onto the $y$ coordinate. The 
geometry of these conditions is shown in Fig.~\ref{fig:covrel}. For the precise 
statement of a general theorem concerning, the method of covering we refer to 
\cite{Zgliczynski1997}. 
  
\begin{figure}[htbp]
 \centerline{\includegraphics[width=\textwidth]{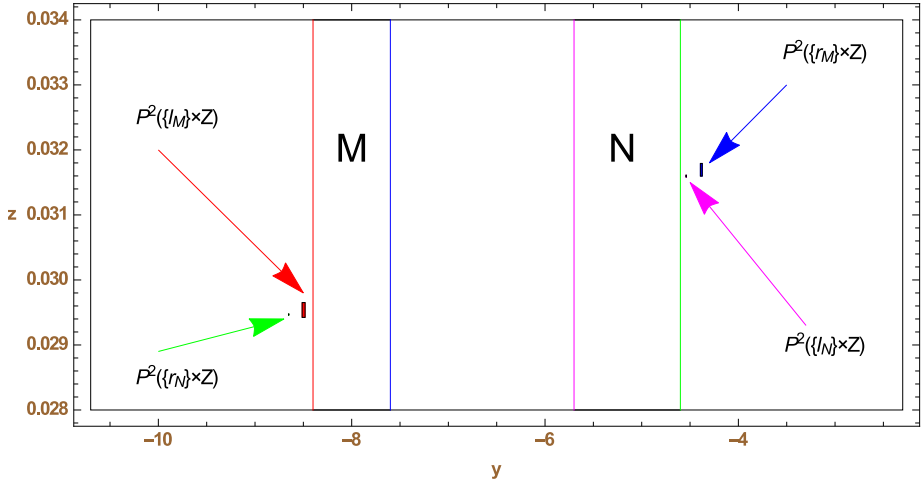}}
 \caption{The sets $M$ and $N$ and rigorous enclosures of the images of their 
exit edges --- see (\ref{eq:covrel}).\label{fig:chaos}}
\end{figure}

The conditions (\ref{eq:covrel}) have been verified in direct computation. We 
did not need to subdivide any of the four edges of $N$ and $M$ that appear in 
(\ref{eq:covrel}). Rigorous bounds on $P^2(\{l_M\}\times Z)$, $P^2(\{r_M\}\times 
Z)$, $P^2(\{l_N\}\times Z)$ and $P^2(\{r_N\}\times Z)$, returned by our 
routine, are shown in Fig.~\ref{fig:chaos}.
  
\textbf{Hyperbolicity and full conjugacy}. Uniform hyperbolicity of $\mathcal 
H$ is proved by means of the cone condition introduced in 
\cite{KokubuWilczakZgliczynski2007}. Here we use our algorithm for integration 
of variational equations. Derivatives with respect to initial conditions are 
necessary for computation of the derivative of Poincar\'e map $P^2$. Let $Q$ be 
a diagonal matrix $Q=\mathrm{Diag}(\lambda,\mu)$ with arbitrary coefficients 
satisfying $\lambda>0$ and $\mu<0$. It has been shown \cite{Wilczak2010} that if 
for all $(y,z)\in N\cup M$ the matrix
 \begin{equation}\label{eq:coneCondition}
  DP^2(y,z)^T\cdot Q \cdot DP^2(y,z) -Q
 \end{equation}
is positive definite, then the maximal invariant set for $P^2$ in $N\cup M$ 
is uniformly hyperbolic. In our computations we used $\lambda=1$ and 
$\mu=-1000$. 

 We uniformly subdivided both sets $N$ and $M$ onto $48$ and $32$ equal pieces, 
respectively (only $y$ coordinate was subdivided). Then, each rectangle was 
submitted to our routine that integrates the first order variational equations 
and computes derivative of the Poincar\'e map $P^2$. Given a rigorous bound of 
the derivative, we checked successfully the condition (\ref{eq:coneCondition}). 
Note that in the case of $2\times 2$ matrix it is easy to check positive 
definiteness of a matrix by the Sylvester criterion. 
\qed

\subsection{Comparison of time of computation.}

In the section, we discuss how the CPU-time needed for verification of the
uniform hyperbolicity in Theorem~\ref{thm:rosslerMainResult} depends on the 
choice of the algorithm used to integrate variational equations. To this end, we 
did the following numerical experiment. For  fixed parameters
\begin{itemize}
 \item $m$ --- the order of numerical method,
 \item $tol$ --- truncation error per one step of the numerical method,
 \item $Alg$ --- the algorithm used to integrate variational equations (\lo or \ho algorithm)
\end{itemize}
we compute the following three numbers
\begin{itemize}
 \item $g_N(m,tol)$, $g_M(m,tol)$ --- minimal natural numbers, such that using 
algorithm $Alg$, the method of order $m$ with the tolerance $tol$ we were able 
to check the cone condition (\ref{eq:coneCondition}) subdividing uniformly the 
sets $N$, $M$ onto $g_N$ and $g_M$ parts, respectively,
 \item $t(Alg)$ --- CPU time of checking the cone condition on both sets $N$ and $M$ with the algorithm and parameters as above.
\end{itemize}

\renewcommand{\arraystretch}{1.2}
\begin{table}[htbp]
\centering

\begin{tabular}{c|c||r|r|r||r|r|r||r}
\hline\hline
 &  & \multicolumn{3}{|c||}{\lo} & \multicolumn{3}{|c||}{\ho}\\
\cline{3-8}
$m$& $tol$& $g_M$ & $g_N$ & $t(\mathrm{LO})$ & $g_M$ & $g_N$ & $t(\mathrm{HO})$ 
& $\displaystyle\frac{t(\mathrm{LO})}{t(\mathrm{HO})}$\\ 
\hline
$10$ & $10^{-10}$ & $39$ & $33$ & $0.90$&$48$ & $32$ & $0.82$& $1.10$\\$10$ & $10^{-12}$ & $38$ & $31$ & $1.29$&$37$ & $31$ & $1.02$& $1.26$\\$10$ & $10^{-14}$ & $25$ & $31$ & $1.69$&$23$ & $30$ & $1.20$& $1.41$\\$10$ & $10^{-16}$ & $25$ & $28$ & $2.25$&$24$ & $25$ & $1.67$& $1.35$\\\hline
$14$ & $10^{-10}$ & $45$ & $52$ & $1.13$&$48$ & $48$ & $0.84$& $1.34$\\$14$ & $10^{-12}$ & $42$ & $39$ & $1.22$&$41$ & $36$ & $0.90$& $1.35$\\$14$ & $10^{-14}$ & $47$ & $33$ & $1.65$&$49$ & $33$ & $1.21$& $1.36$\\$14$ & $10^{-16}$ & $36$ & $32$ & $1.70$&$36$ & $31$ & $1.38$& $1.23$\\\hline
$18$ & $10^{-10}$ & $63$ & $77$ & $1.67$&$62$ & $56$ & $1.20$& $1.40$\\$18$ & $10^{-12}$ & $48$ & $56$ & $1.41$&$63$ & $49$ & $1.26$& $1.12$\\$18$ & $10^{-14}$ & $44$ & $43$ & $1.76$&$47$ & $38$ & $1.18$& $1.50$\\$18$ & $10^{-16}$ & $40$ & $37$ & $1.91$&$54$ & $34$ & $1.44$& $1.33$\\\hline
$22$ & $10^{-10}$ & $151$ & $95$ & $3.36$&$101$ & $67$ & $1.93$& $1.74$\\$22$ & $10^{-12}$ & $78$ & $78$ & $2.44$&$59$ & $58$ & $1.81$& $1.35$\\$22$ & $10^{-14}$ & $52$ & $61$ & $2.05$&$61$ & $49$ & $1.56$& $1.32$\\$22$ & $10^{-16}$ & $45$ & $41$ & $1.80$&$47$ & $47$ & $1.53$& $1.17$\\\hline
\hline
\end{tabular}
\caption{Comparison of \lo and \ho algorithms.\label{tab:rosslerData}}
\end{table}

Let us emphasize that the vector field of the R\"ossler system 
(\ref{eq:rosslerSystem}) contains only one nonlinear term. Hence, we have 
$c_f=1$, and this is almost the worst linear case for the \ho method when the 
complexity is dominated by the matrix operations and the expected time 
savings from the \ho method are smaller --- see analysis in 
Section~\ref{sec:complexity}.  

In Table~\ref{tab:rosslerData}, we present results from this experiment. We see 
that in each case the \ho algorithm is faster than the \lo algorithm. Higher 
computational complexity of the \ho algorithm is compensated by significantly 
smaller truncation error. Therefore, a routine that predicts the time step (the 
same routine was used in both cases) returns larger time steps for the \ho 
algorithm, and in consequence, the total computing time is smaller. We also 
notice that in some cases decreasing the tolerance increases the number of 
subdivisions $g_N$ and $g_M$ needed to check the cone condition --- see for 
instance 
the row with $m=14$ and $tol=10^{-14}$. This is a consequence of many heuristics 
made in the implementation (for instance reorganization of doubleton 
representation after reaching some threshold values). These heuristics make the 
algorithm discontinuous with respect to parameters. Moreover, decreasing the 
tolerance increases number of time steps needed to compute a full trajectory. 
This may result in weaker control of unavoidable wrapping effect.

\section{Conclusions.}
Since the \lo algorithm appeared \cite{Zgliczynski2002} it has been proved to be 
very useful in rigorous analysis of ODEs. In this paper, we proposed an 
efficient 
alternative for this algorithm and we provided free implementation of both \lo 
and \ho algorithms available as a module of the CAPD library \cite{CAPD}. 
Numerical tests show that the \ho algorithm is slightly faster than the widely 
used \lo algorithm. We have shown that the \ho algorithm may be faster in 
practical applications. This is not very important when the total time of 
computation is counted in seconds, as we have seen in 
Section~\ref{sec:applications}. Any progress matters, however, if a problem 
requires hundreds or thousands CPU hours: for example verification of the 
existence of an uniformly hyperbolic attractor of the Smale-Williams type 
\cite{Wilczak2010} or the coexistence of chaos and hyperchaos in the 4D 
R\"ossler system 
\cite{BarrioMartinezSerranoWilczak2015,WilczakSerranoBarrio2016}. In the 
computation reported in 
\cite{BarrioMartinezSerranoWilczak2015,WilczakSerranoBarrio2016}, the proposed 
\ho algorithm has been used.

In \cite{WilczakZgliczynski2011}, an algorithm for integration of higher order 
variational equations is presented. Ideas from Section~\ref{sec:algorithm} can 
be directly used to design and implement an algorithm, let us call it $\mathcal 
C^r$-HO, with the $\mathcal C^r$-Lohner method as a predictor step. This 
requires encoding rather than theoretical effort and, we hope, this 
implementation will be available soon as part of the CAPD library \cite{CAPD}.

\end{document}